\documentclass[final]{siamltex}

\usepackage{amscd,amsmath,amssymb,epsfig}
\usepackage{tikz}

\newtheorem{example}[theorem]{Example}
\newtheorem{remark}[theorem]{Remark}

\title{Max-plus objects to study the complexity of graphs}

\author{Cristiano Bocci\thanks{Dipartimento di Scienze Matematiche e Informatiche ``R. Magari'', Universit\`a di Siena, Pian dei Mantellini 44, 53100 Siena, Italy ({\tt cristiano.bocci@unisi.it})}
        \and Luca Chiantini\thanks{Dipartimento di Scienze Matematiche e Informatiche ``R. Magari'', Universit\`a di Siena, Pian dei Mantellini 44, 53100 Siena, Italy ({\tt luca.chiantini@unisi.it})} \and Fabio
Rapallo\thanks{Dipartimento di Scienze e Tecnologie Avanzate,
Universit\`a del Piemonte Orientale, Viale Teresa Michel 11, 15121
Alessandria, Italy ({\tt fabio.rapallo@mfn.unipmn.it})}}

\begin{document}

\maketitle

\begin{abstract}
Given an undirected graph $G$, we define a new object $H_G$,
called the mp-chart of $G$, in the max-plus algebra. We use it,
together with the max-plus permanent, to describe the complexity
of graphs. We show how to compute the mean and the variance of
$H_G$ in terms of the adjacency matrix of $G$ and we give a
central limit theorem for $H_G$. Finally, we show that the
mp-chart is easily tractable also for the complement graph.
\end{abstract}



\pagestyle{myheadings} \thispagestyle{plain} \markboth{C. BOCCI,
L. CHIANTINI AND F. RAPALLO}{MAX-PLUS OBJECTS FOR GRAPHS}

\section{Introduction}

The work presented in this paper has been inspired by the need of
simple and actual techniques to measure the complexity of a graph,
especially in the case of sparse graphs. This problem arises in
several fields of applications, from Computer Science to
Economics, from Biology to Social Sciences. As general references
for the graph theory, we mention in particular the books
\cite{bollobas:98}, \cite{bollobas:01}, \cite{biggs:93} and
\cite{diestel:05}, where the reader can find the main mathematical
achievements in the theory. For a general survey on recent
applications of graph theory see, for instance,
\cite{albert|barabasi:02}. The reader interested in some more
technical papers can refer to \cite{keller:07} for applications to
Economics, to \cite{chakraboti|etal:10} for applications to
Econophysics, to \cite{mason|verwoerd:07} for applications to
Biology, and to \cite{hwang|etal:08} for applications to Molecular
Biology. In such papers, sparse graphs play a prominent role.

Although an undirected graph $G=(V,E)$ is a rather simple
structure, consisting of a set $V$ of $N$ vertices and a set of
edges $E \subset V \times V$, in graph theory there are several
different approaches, depending on the specific application we are
looking for. In particular, a graph can be fixed or random,
depending on whether the elements in $E$ are random or not.
Moreover, many efforts have been done to analyze dynamical graphs,
where the vertex set $V$ and/or the edge set $E$ vary with time,
see e.g. \cite{durrett:07}.

Here, we restrict our analysis to fixed graphs. In this framework,
there are interesting developments in the area of Combinatorics,
about the study of the properties of $0-1$ matrices. These
matrices naturally arise in the framework of graphs, as the
adjacency matrix $A_G$ of the fixed graph $G$. Some recent
developments in this direction, with applications to graph theory,
are described in \cite{barvinok:08}, \cite{barvinok:10}, and
\cite{barvinok|hartigan:10}.

In the present paper, we investigate some questions about
undirected graph, in order to study and describe their structure,
with special attention to sparse graphs. Our work is related to
the matching problem. As a preliminary remark, we argue that, for
sparse graphs, the classical descriptors of the complexity, such
as the degree distribution and the permanent of the adjacency
matrix do not give actual information. Thus, we use the max-plus
arithmetic and the corresponding expression of the permanent, and
we show that this object is more suitable for describing of the
complexity for sparse graph. The use of the max-plus arithmetic
naturally leads to the definition of a more complete index of the
structure of a graph, and we define a vector called the {\it
mp-chart of the graph}. This vector is nothing else but a
probability distribution, and we show that it converges to a
Normal distribution through a combinatorial Central Limit Theorem.
Several examples on small- and medium-sized graphs are given to
show that our definitions are easy to apply and provide practical
information about the complexity of the structure of the graph
under study. All the computations have been carried out with
Maple, see \cite{maple}, and {\tt R}, see \cite{rproject}. All the
simulations come from simple {\tt R} routines, without any
additional package.

This paper is only concerned with undirected fixed graphs.
Nevertheless, the same strategy can be applied to other
situations, such as bipartite graphs, undirected unfixed graphs,
random graphs, and so on.

The paper is organized as follows. In Section \ref{mppermsect} we
define the max-plus permanent of an adjacency matrix (i.e., the
permanent under the max-plus arithmetic), we state its main
properties, we study its connections with the classical permanent,
and we discuss some simple examples. In Section \ref{mpchartsect}
we define a new object associated to a graph, and we call it the
mp-chart of the graph. We compute its mean and variance, and we
show that, under suitable conditions, it converges to a Normal
distribution as the size of the graph goes to infinity. Some
simple simulations show that the convergence is quite good also
for small values of the size. In Section \ref{complsect}, we show
how the mp-chart of a graph is related to the mp-chart of the
complement graph. Finally, Section \ref{futuresect} is devoted to
suggest some future directions of this research.

\section{The max-plus permanent} \label{mppermsect}

Let $G=(V,E)$ be an undirected graph with $N$ vertices. Let $A_G$
be the $N \times N$ adjacency matrix of $G$, defined by
$(A_G)_{i,j}=1$ if $(i,j) \in E$ and $0$ otherwise. In the classic
definition of undirected graph, the matrix $A_G$ is symmetric and
with zero diagonal entries, as we do not consider loops.

As mentioned in the Introduction, the study of the complexity of a
given graph is one of the most relevant problems about graphs in
Applied Probability. This analysis can be performed through the
distribution of the degrees (i.e., the number of edges involving
each vertex) and through the permanent (or the determinant) of the
adjacency matrix $A_G$.

The determinant of $A_G$ is
\begin{equation*}
\mathrm{det}(A_G)= \sum_{\pi} (-1)^{|\pi|} \prod_{i=1}^N
(A_G)_{i,\pi(i)} \, ,
\end{equation*}
where the sum is taken over all the permutations $\pi$ of $\{1,
\ldots, N\}$ and $|\pi|$ denotes the parity of $\pi$. The
permanent of $A_G$ is
\begin{equation*}
\mathrm{perm}(A_G)= \sum_{\pi} \prod_{i=1}^N (A_G)_{i,\pi(i)} \, .
\end{equation*}

The use of permanent to describe the complexity of a graph is
justified by the following well-known property.

\begin{proposition}
The permanent of $A_G$ is the number of bijections $\phi:V
\rightarrow V$ compatible with $E$, i.e. such that $(v,\phi(v))
\in E$ for all $v \in V$.
\end{proposition}

In fact, $\mathrm{perm}(A_G)$ is the number of permutations $\pi$
with $A_{1,\pi(1)}= \ldots = A_{N,\pi(N)}=1$ and the permutation
$\pi$ is just the bijection $\phi$ in the proposition.

However, the analysis based on the degree distribution and the
permanent is not adequate for sparse graphs. In fact, it is enough
to have an isolated vertex to produce a null permanent.
Nevertheless, it is interesting to study the structure of a sparse
graph.

To overcome this difficulty, we make use of the tropicalization of
the permanent. In the classical settings, Tropical Arithmetic is
defined through the operations:
\begin{equation*}
x \oplus y = \min\{x,y\} \qquad x \otimes y = x + y
\end{equation*}

But, with Tropical Arithmetic, the determinant (or permanent) of
an adjacency matrix is always $0$, because of the nullity of the
main diagonal of $A_G$.

Thus, we use the max-plus algebra, with operations:
\begin{equation*}
x \oplus y = \max\{x,y\} \qquad x \otimes y = x + y
\end{equation*}

Consequently, the explicit expression of the max-plus permanent is
\begin{equation} \label{mppermdef}
\mathrm{perm_{mp}}(A_G):= \bigoplus_{\pi} \left(
\bigotimes_{i=1}^N (A_G)_{i,\pi(i)} \right) = \max_{\pi}
\sum_{i=1}^N (A_G)_{i,\pi(i)}
\end{equation}

The max-plus permanent is the maximum over $N!$ terms. Each of the
$N!$ terms is the sum of $N$ terms in $\{0,1\}$. Thus, the
max-plus permanent $\mathrm{perm_{mp}}(A_G)$ is zero if and only
if the matrix $A_G$ is the null matrix. On the opposite side, the
maximum allowed value of the max-plus permanent is $N$.

The use of the max-plus permanent to analyze sparse graphs has a
first reason in the following property.
\begin{lemma} \label{mp-perm-inter}
The following relation holds:
\begin{equation}
\mathrm{perm_{mp}}(A_G)=N \Longleftrightarrow \mathrm{perm}(A_G)>0
\, .
\end{equation}
\end{lemma}
\begin{proof}
$\mathrm{perm_{mp}}(A_G)=N$ if and only if there exists a
permutation $\pi$ such that $\sum_{i=1}^N (A_G)_{i,\pi(i)}=N$.
This happens if and only if there exists $\pi$ such that
$(A_G)_{i,\pi(i)}=1$ for $i=1, \dots, N$, i.e. if and only if $
\mathrm{perm}(A_G)>0$.
\end{proof}

\begin{remark}
Notice that, from Lemma \ref{mp-perm-inter} and from the previous
discussion, it follows that the max-plus permanent is able to
discriminate among graphs with standard permanent equal to zero.
\end{remark}

Moreover, we explicitly write the following consistency property,
whose simple proof is straightforward.

\begin{lemma}
Let $G$ and $H$ be two graphs on two disjoint sets of vertices.
Then,
\begin{equation}\label{disjoint}
\mathrm{perm_{mp}}(A_{G\cup
H})=\mathrm{perm_{mp}}(A_G)+\mathrm{perm_{mp}}(A_H) \, .
\end{equation}
\end{lemma}

The max-plus permanent has interesting connections with the
subgraphs. Let $G'=(V',E')$ a graph. If $V' \subset V$ and $E'
\subset E$, then $G'$ is a subgraph of $G=(V,E)$. A subgraph
$G'=(V',E')$ is the subgraph induced by $V'$ if $E'$ contains all
the edges in $E$ involving the vertices in $V'$. In order to
analyze the max-plus permanent, in view of Equation
\eqref{mppermdef}, we introduce here the notion of $t$-term, which
is strictly related to the subgraphs of $G$. Such connections will
be studied later in this section.

\begin{definition}
Given a graph $G$ with adjacency matrix $A_G$, a $t$-term is a
sequence of indices $(i_1,j_1) \cdots (i_t,j_t)$ with
\begin{itemize}
\item $1 \leq i_1 < \ldots < i_t \leq N$;

\item the $j_k$'s, with $1 \leq j_k \leq N$ are all distinct;

\item $(A_G)_{i_k,j_k}=1$ for all $k$.
\end{itemize}
For a $t$-term $P$, we denote $I(P)=\{ i_1, \dots, i_t\}$ and
$J(P)=\{ j_1, \dots, j_t\}$.
\end{definition}

Roughly speaking, a $t$-term corresponds to a sequence of
positions of $t$ ones in the permutations. A straightforward
consequence is the following statement.

\begin{proposition}
The max-plus permanent of $G$ is $q$ if and only if there exists a
$q$-term and there are no $t$-term with $t > q$.
\end{proposition}

\begin{proposition} \label{prop-subgraphs}
Let $t$ the maximum integer such that there exists a $t$-term,
then there exists a $t$-term $P$ such that $I(P)=J(P)$.
\end{proposition}

\begin{proof}
We prove the statement by induction on $t$. If $t=2$ there is
nothing to prove since if the $2$-term is given by
$(i_1,j_1)(i_2,j_2)$ it is enough to consider  the $2$-term
$(i_1,j_1)(j_1,i_1)$.

Consider now a $t$-term $Q$ and suppose there exists a $k$ such
that $i_k$ is in $I(Q) \setminus J(Q)$. First of all we notice
that $j_k$ must be in $I(Q)$. If not, we can add the element
$(j_k,i_k)$ to $Q$ obtaining a $(t+1)$-term which is a
contradiction, since $Q$ is maximal. Hence, since $j_k \in I(Q)$
then there exists a $s$ such that $j_k=i_s \in I(Q)$. Then, we
substitute $(i_s,j_s)$ with $(j_k,i_k)$ in our $t$-term and we
obtain a new $t$-term of the form $(i_k,j_k)(j_k,i_k)Q'$ where
$Q'$ is a $(t-2)$-term. This term $Q'$ arises from the sub-matrix
$A'$ of $A(G)$ where we remove rows and columns $i_k$ and $j_k$.
Hence $Q'$ is a maximal $(t-2)$-term for $A'$. By induction, the
proof follows.
\end{proof}

\begin{remark}
If $P$ is a $3$-term, then we must have $I(P)=J(P)$. In fact, if
$P$ is $(i_1,j_1)(i_2,j_2)(i_3,j_3)$ with $I(P)\not=J(P)$ then, by
the previous proposition, we obtain a new $3-$term
$(i_k,j_k)(j_k,i_k)(i_3,j_3)$. Then it would be possible to extend
it to $(i_k,j_k),(j_k,i_k)(i_3,j_3)(j_3,i_3)$ against the
maximality of $P$.
\end{remark}

\begin{remark}
In view of Proposition \ref{prop-subgraphs}, the max-plus
permanent is the cardinality of the largest subset of $V$ with a
bijection compatible with $E$. This is another way to see that the
max-plus permanent is able to detect the complexity of the graphs
with null classical permanent.
\end{remark}

Denote by $\ell_G$ the number of edges of a graph $G$. Among the
subgraphs of Proposition \ref{prop-subgraphs}, we are mainly
interested in the ones with a minimal number of edges. These
subgraphs are maximal in term of $\mathrm{perm_{mp}}{(A_{G'})}$,
but minimal in term of $\ell_{G'}$. We made this more precise by
the following definition.

\begin{definition}
An  {\it mp-maximal subgraph} $G'$ of a graph $G$, is a  subgraph
of $G$ with $q=\mathrm{perm_{mp}}(A_{G})$ vertices,
\begin{equation}\label{maxmin}
\mathrm{perm_{mp}}(A_{G'})=\mathrm{perm_{mp}}(A_G)
\end{equation}
and for all other subgraph $G''$ of $G$ satisfying (\ref{maxmin})
one has $\ell_{G'} \leq \ell_{G''}$.
\end{definition}

The rest of this section is devoted to the discussion of some
examples and some useful remarks. In order to understand the
definitions introduced above, we start with some small graphs.

\begin{example} \label{firstex}
Let us analyze the three graphs on $4$ vertices drawn in Figure
\ref{firstfig}. Their adjacency matrices are respectively
\begin{equation*}
A_{G_1}=\begin{pmatrix}
0 & 1 & 0 & 1\\
1 & 0 & 1 & 0\\
0 & 1 & 0 & 1\\
1 & 0 & 1 & 0
\end{pmatrix}
\qquad A_{G_2}=\begin{pmatrix}
0 & 1 & 0 & 0\\
1 & 0 & 0 & 0\\
0 & 0 & 0 & 1\\
0 & 0 & 1 & 0
\end{pmatrix}
\qquad A_{G_3}=\begin{pmatrix}
0 & 1 & 1 & 0\\
1 & 0 & 1 & 0\\
1 & 1 & 0 & 0\\
0 & 0 & 0 & 0
\end{pmatrix}
\end{equation*}

In the first graph, all vertices are connected and
$\mathrm{perm_{mp}}(A_{G_1})=4$. However, this is not the minimal
way to obtain a max-plus permanent equal to $4$. In fact, it is
easy to check that $\mathrm{perm_{mp}}(A_{G_2})=4$. Thus the graph
$G_2$ represents a mp-maximal subgraph for the graph $G_1$, but it
is not the only one. If we look now at the graph $G_3$, we notice
that there is a cycle of length $3$ and an isolated vertex. In
such case, we have $\mathrm{perm_{mp}}(A_{G_3})=3$, and there is
only one mp-maximal subgraph.
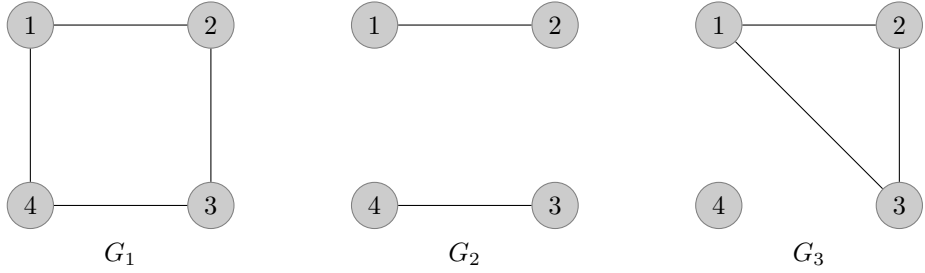
\begin{figure}
\begin{center}
\begin{tabular}{ccccc}
\begin{tikzpicture}
  [scale=.6,auto=left,every node/.style={circle,draw=black!50,fill=black!20}]
  \node (n1) at (0,4) {1};
  \node (n2) at (4,4) {2};
  \node (n3) at (4,0) {3};
  \node (n4) at (0,0) {4};

  \foreach \from/\to in {n1/n2,n2/n3,n3/n4,n4/n1}
    \draw (\from) -- (\to);
\end{tikzpicture}
& \qquad \qquad &
\begin{tikzpicture}
  [scale=.6,auto=left,every node/.style={circle,draw=black!50,fill=black!20}]
  \node (n1) at (0,4) {1};
  \node (n2) at (4,4) {2};
  \node (n3) at (4,0) {3};
  \node (n4) at (0,0) {4};

  \foreach \from/\to in {n1/n2,n3/n4}
    \draw (\from) -- (\to);
\end{tikzpicture}
& \qquad \qquad &
\begin{tikzpicture}
  [scale=.6,auto=left,every node/.style={circle,draw=black!50,fill=black!20}]
  \node (n1) at (0,4) {1};
  \node (n2) at (4,4) {2};
  \node (n3) at (4,0) {3};
  \node (n4) at (0,0) {4};

  \foreach \from/\to in {n1/n2,n1/n3,n2/n3}
    \draw (\from) -- (\to);
\end{tikzpicture} \\
$G_1$ &  & $G_2$  & & $G_3$
\end{tabular}
\end{center}
\caption{The graphs $G_1$, $G_2$ and $G_3$ for Example
\ref{firstex}.} \label{firstfig}
\end{figure}
\end{example}

\begin{example} \label{oppositeex}
To illustrate the behavior of the max-plus permanent and of the
mp-maximal subgraphs, we analyze two opposite examples.  with the
same length. The two graphs are drawn in Figure
\ref{esempiopposti}. The graph $G_1$ on the left is the union of a
tree and two isolated vertices, with
$\mathrm{perm_{mp}}(A_{G_1})=2$ and $4$ maximal subgraphs with two
vertices and one edge each. On the opposite side, the graph $G_2$
has a perfect matching, $\mathrm{perm_{mp}}(A_{G_2})=6$ and there
is only $1$ maximal subgraph, i.e., the graph $G_2$ itself.

\begin{figure}
\begin{center}
\begin{tabular}{ccc}
\begin{tikzpicture}
  [scale=.6,auto=left,every node/.style={circle,draw=black!50,fill=black!20}]
  \node (n1) at (0,4) {1};
  \node (n2) at (3,4) {2};
  \node (n3) at (6,4) {3};
  \node (n4) at (0,0) {4};
  \node (n5) at (3,0) {5};
  \node (n6) at (6,0) {6};

  \foreach \from/\to in {n2/n4,n2/n5,n2/n6}
    \draw (\from) -- (\to);
\end{tikzpicture}
& \qquad \qquad &
\begin{tikzpicture}
  [scale=.6,auto=left,every node/.style={circle,draw=black!50,fill=black!20}]
  \node (n1) at (0,4) {1};
  \node (n2) at (3,4) {2};
  \node (n3) at (6,4) {3};
  \node (n4) at (0,0) {4};
  \node (n5) at (3,0) {5};
  \node (n6) at (6,0) {6};

  \foreach \from/\to in {n1/n4,n2/n5,n3/n6}
    \draw (\from) -- (\to);
\end{tikzpicture} \\
$G_1$ &  & $G_2$
\end{tabular}
\end{center}
\caption{The two graphs $G_1$ and $G_2$ for Example
\ref{oppositeex}.} \label{esempiopposti}
\end{figure}
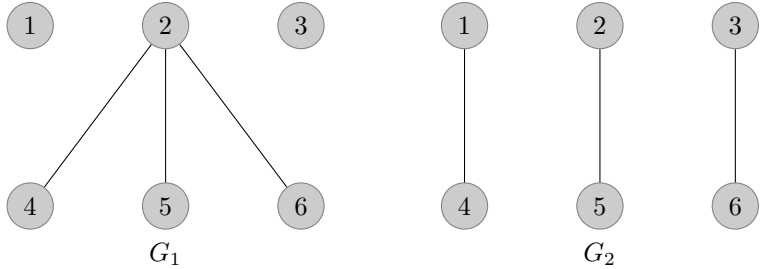
\end{example}

%\begin{Rmk}
%Let $t=\mathrm{perm_{mp}}(A_{G})$, then each mp-maximal subgraph
%corresponds to a $t$-term as built in the proof of Proposition
%\ref{prop-subgraphs}.
%\end{Rmk}

\begin{proposition}
Two mp-maximal subgraphs are not disjoint.
\end{proposition}
\begin{proof}
Consider a graph $G$ such that $\mathrm{perm_{mp}}(A_G)=t$ and let
$G'=(V',E')$ and $G''=(V'',E'')$ be two mp-maximal subgraphs of
$G$ with $t$ vertices each. Suppose that $V'$ and $V''$ are
disjoint. Then, by Formula (\ref{disjoint}), the adjacency matrix
of $G'\cup G''$ has max-plus permanent $2t$. Then
$\mathrm{perm_{mp}}(A_G)\geq 2t$ which is a contradiction.
\end{proof}

\begin{remark}
In the max-plus arithmetic, the definition of determinant is not
unique, see \cite{gaubert|meunier:10}. Therefore, one has to
define the positive and negative determinant. In particular, the
positive max-plus determinant is the maximum of the sums $\sum_{i}
A_{i, \pi(i)}$  over all even permutations $\pi$. The negative
determinant is defined by taking the odd permutations instead of
the even ones. This issue is another reason to use the permanent
instead of the determinant in the max-plus environment.
\end{remark}

%%%%%%%%%%%%%%%%%%%%%%%%%%%%%%%%%%%%%%%%%%%%%%%%%%%%%%%%%%%

\section{The mp-chart of a graph} \label{mpchartsect}

The information about a graph is not contained only in the
max-plus permanent, but in the whole distribution of the $N!$
terms $\sum_{i=1}^N (A_G)_{i,\pi(i)}$. Thus, in this section we
define the mp-chart of a graph as the distribution of the $N!$
terms above, and we prove that this distribution converges to a
Gaussian distribution through the Hoeffding's combinatorial
central limit theorem, see \cite{hoeffding:51}. We also show that
the mean and the variance of that distribution can be computed
easily from the adjacency matrix.

\begin{definition}
Let $G$ be a graph and $A_G$ its adjacency matrix. Let $h_G(k)$ be
the number of permutations $\pi$ such that $\sum_{i=1}^N
(A_G)_{i,\pi(i)}=k$. We call the $(N+1)$-dimensional integer
vector $H_G=(h_G(0), \dots, h_G(N))$ the {\it mp-chart of the
graph} $G$.
\end{definition}

This object captures many features of the graph and has some
relevant theoretical properties. To understand the meaning of
$H_G$, notice that  $h_G(k)$ is just the number of permutations
$\pi$ such that the sequence $(A_G)_{1,\pi(1)}, \ldots ,
(A_G)_{N,\pi(N)}$ contains a $k$-term but not a $(k+1)$-term. This
gives precisely the meaning and the usefulness of the notion of
random permutation in that context.

\begin{example} We use here a very simple scheme inspired by Econophysics, see
\cite{silver|slud|takamoto:02} and \cite{garibaldi|scalas:10}.
Consider a population with $N$ agents, each possessing one good.
The goods can be sent and received only along the edges of a graph
$G$ and each agent can possess only one good. Given a random
permutation $\pi$ of $\{1, \ldots, N\}$, the $i$-th agent can send
its good to $\pi(i)$ if it receive a good from $\pi^{-1}(i)$. The
quantity $\sum_{i=1}^N (A_G)_{i, \pi(i)}$ is exactly the number of
agents involved in this process. Of course, similar examples can
be adapted to many other sciences.
\end{example}

We start the analysis of the mp-chart with the study of the mean
$\mathbb{E}(H_G)$ and the variance $\mathbb{V}(H_G)$. Although
these computations could be carried out applying Theorem 2 in
\cite{hoeffding:51}, it is useful to state explicitly the proof
for adjacency matrices.

Given an adjacency matrix $A=(A_G)$, note that, for a random
permutation $\pi$, the object
\begin{equation}
S_N(\pi) = A_{1,\pi(1)} + \ldots + A_{N,\pi(N)}
\end{equation}
is the sum of $N$ binary random variables, and $A_{i,\pi(i)}$ is
the addendum chosen in the $i$-th row of the adjacency matrix.

In order to analyze the mean of $H_G$, defined by
\begin{equation}
\mathbb{E}(H_G)= \frac {\sum_{t=0}^N kh_G(k) } {N!} \, ,
\end{equation}
it is convenient to adopt an inductive approach.

\begin{theorem}\label{Media}
The mean of the mp-chart $H_G$ is
\begin{equation} \label{mean}
\mathbb{E}(H_G)=\frac{2 \ell_G}{N} \, ,
\end{equation}
where $\ell_G$ is the number of edges in the graph $G$.
\end{theorem}
\begin{proof}
Clearly, if $A_G$ is the null matrix, then $\mathbb{E}(H_G)=0$.
Suppose that the formula \eqref{mean} holds true for $\ell_G-1$.
By direct inspection, adding one edge has the following
consequences. Among the $N!$ terms $S_N(\pi)$:
%\begin{figure}
%\begin{center}
%$\left(\begin{array}{cc|ccccc}
%0 & 1 & & \cdots & & \cdots & \\
%1 & 0 & & \cdots & & \cdots &  \\ \hline
%  &   & & & & & \\
%\vdots  &  \vdots & & & & & \\
%& & & & & \\
%\vdots & \vdots& & & &\\
%&
%\end{array}\right)
%$ \caption{Adding the edge $(1,2)$ means adding two ones in the
%position $(1,2)$ and $(2,1)$ of the adjacency matrix.}
%\label{add-edge}
%\end{center}
%\end{figure}
\begin{itemize}
\item $(N-2)!$ of them increase by $2$;

\item $2(N-2)(N-2)!$ of them increase by $1$;

\item the remaining $((N-2)^2+N-1)(N-2)!$ do not change.
\end{itemize}
Thus,
\begin{equation*}
\mathbb{E}(H_G)=\frac{2 (\ell_G-1)}{N} + \frac
{2(N-2)!+2(N-2)(N-2)!} {N!} = \frac{2 (\ell_G-1)}{N} + \frac 2 N =
\frac{2\ell_G} N \, .
\end{equation*}
\end{proof}

\begin{example}
Given a complete graph $G$, its adjacency matrix has $0$ on the
diagonal and $1$ elsewhere. The graph has ${N(N-1)}/{2}$ edges.
Hence $\mathbb{E}(H_G)=N-1$ which is the maximum allowed.
\end{example}

Notice that the mean $\mathbb{E}(H_G)$ depends only in the number
of edges of $G$, whatever they are collocated, that is,
$\mathbb{E}(H_G)$ does not take into account the topology of the
graph. On the other hand, the variance $\mathbb{V}(H_G)$ depends
on the position of the edges.

\begin{theorem}\label{Varianza}
The variance of the mp-chart $H_G$ is
\begin{equation} \label{variance}
\mathbb{V}(H_G)= \sum_{i=1}^N \frac {d_i(N-d_i)} {N^2} + \sum_{
\substack{i,j=1 \\ i \ne j} }^N \frac {d_id_j - NT_{i,j}}
{N^2(N-1)} \, ,
\end{equation}
where $d_1, \ldots , d_N$ are the degrees of the vertices and
$T_{i,j}=<r_i,r_j>$ is the scalar product of the $i$-th and the
$j$-th row of $A_G$.
\end{theorem}

\begin{proof}
By direct computation, the formula \eqref{variance} holds for $N
\leq 2$.

To prove the validity of Eq. \eqref{variance} for $N \geq 3$, it
is enough to compute the covariances
\begin{eqnarray*}
\mathrm{Cov}(A_{i,\pi(i)},A_{j,\pi(j)}) &=&
\mathbb{E}(A_{i,\pi(i)} A_{j,\pi(j)}) - d_id_j/N^2 = \\
&=& \mathbb{P}(A_{i,\pi(i)}=1, A_{j,\pi(j)}=1) - d_id_j/N^2  \, .
\end{eqnarray*}
Without loss of generality we can fix $i=1$ and $j=2$ and we write
for brevity $A_1$ for $A_{1,\pi(1)}$ and $A_2$ for $A_{2,\pi(2)}$.
Moreover, we suppose that $d_1$ and $d_2$ are both non zero. (If
$d_1=0$ or $d_2=0$, then trivially $\mathrm{Cov}(A_1,A_2)=0)$.

We divide the computation in two cases, and to help the reader we
have sketched the two cases in Figure \ref{fig-dim}.

\begin{figure}
\begin{center}
\begin{tabular}{cc}
$\left(\begin{array}{cc|ccccc}
0 & 0 & & \cdots & & \cdots & \\
0 & 0 & & \cdots & & \cdots &  \\ \hline
  &   & & & & & \\
\vdots  &  \vdots & & & & & \\
& & & & & \\
\vdots & \vdots& & & &\\
&
\end{array}\right)
$ & $\left(\begin{array}{cc|ccccc}
0 & 1 & & \cdots & & \cdots & \\
1 & 0 & & \cdots & & \cdots &  \\ \hline
  &   & & & & & \\
\vdots  &  \vdots & & & & & \\
& & & & & \\
\vdots & \vdots& & & &\\
&
\end{array}\right)
$
\\
Case $(a)$ & Case $(b)$ \end{tabular}
 \caption{The two cases arising in the proof of Theorem \ref{Varianza}.} \label{fig-dim}
\end{center}
\end{figure}

\begin{itemize}
\item Case $(a)$: $(1,2)$ is not an edge of the graph.Then:
\begin{itemize}
\item there are $2(N-1)!$ permutations such that $\pi(1)=1$ or
$\pi(1)=2$. For all these cases, $(A_{1}=1, A_{2}=1)$ is
impossible;

\item there are $2(N-2)(N-2)!$ permutations such that $\pi(1)>2$,
but $\pi(2)=1$ or $\pi(2)=2$. Also in all these cases, $(A_{1}=1,
A_{2}=1)$ is impossible;

\item there are $(N-2)(N-3)(N-2)!$ permutations such that
$\pi(1)>2$ and $\pi(2)>2$, and among these permutations
\begin{equation*}
(T_{1,2} d_2 + (d_1-T_{1,2})d_2)(N-2)! = (d_1d_2 - T_{1,2})(N-2)!
\end{equation*}
are such that $A_{1}A_{2}=1$.
\end{itemize}
Therefore,
\begin{equation*}
\mathrm{Cov}(A_{1},A_{2}) = \frac {d_1d_2 - NT_{1,2}} {N^2(N-1)}
\end{equation*}

\item Case $(b)$: $(1,2)$ is an edge of the graph. Then:
\begin{itemize}
\item there are $(N-1)!$ permutations such that $\pi(1)=1$. In all
such cases, $(A_{1}=1, A_{2}=1)$ is impossible;

\item there are $(N-2)!$ permutations such that $\pi(1)=2$ and
$\pi(2)=1$. For such permutations, $A_{1}=1$ and $A_{2}=1$;

\item there are $(N-2)(N-2)!$ permutations such that $\pi(1)=2$
and $\pi(2)>2$. Among these permutations, $(d_2-1)(N-2)!$ are such
that $A_1A_2=1$.

\item there are $(N-2)(N-2)!$ permutations such that $\pi(1)>2$
and $\pi(2)=1$. Among these permutations, $(d_1-1)(N-2)!$ are such
that $A_1A_2=1$.

\item there are $(N-2)(N-2)!$ permutations such that $\pi(1)>2$
and $\pi(2)=2$. In all such cases, $A_1A_2=0$.

\item there are $(N-2)(N-3)(N-2)!$ permutations such that
$\pi(1)>2$ and $\pi(2)>2$, and among these
\begin{equation*}
\begin{split}
(T_{1,2} (d_2-2) + (d_1-1-T_{1,2})(d_2-1))(N-2)! = \\ = (d_1d_2 -
d_1 -d_2 - T_{1,2} +1)(N-2)!
\end{split}
\end{equation*}
are such that $A_{1}A_{2}=1$.
\end{itemize}
Therefore, adding up all the contributions, we obtain again
\begin{equation*}
\mathrm{Cov}(A_{1},A_{2}) = \frac {d_1d_2 - NT_{1,2}} {N^2(N-1)}
\end{equation*}
\end{itemize}
The formula in Eq. \eqref{variance} is now straightforward.
\end{proof}

\begin{example}
Consider the matrices
\begin{equation*}
A_{G_1}= \begin{pmatrix}
0 & 1 & 0 &0&0\\
1 & 0 & 1 & 0 &0\\
0& 1&0 &0 &0\\
0& 0& 0 &0 &0 \\
0 &0 &0 &0 &0
\end{pmatrix} \qquad A_{G_2}= \begin{pmatrix}
0 & 1 & 0 &0&0 \\
1 & 0 & 0 & 0 &0 \\
0& 0&0 &1 &0 \\
0& 0& 1 &0 &0 \\
0 &0 &0 &0 &0
\end{pmatrix}
\end{equation*}

The graph $G_1$ has two consecutive edges, while the graph $G_2$
has two disjoint edges. An easy computation gives
\begin{equation*}
H_{G_1}=(48, 48, 24,0,0,0)
\end{equation*}
and
\begin{equation*}
H_{G_2}=(53,44,18,4,1,0)
\end{equation*}
with equal means
$\mathbb{E}(H_{G_1})=\mathbb{E}(H_{G_2})={4}/{5}$. On the
contrary, the variances are $\mathbb{V}(H_{G_1})={14}/{25}$ and
$\mathbb{V}(H_{G_2})={19}/{25}$, respectively.
\end{example}

\begin{example}
As a second example, consider the two graphs on the set vertices
$V=\{1 , \ldots, 7 \}$ shown in Figure \ref{fig-eqmv}.
\begin{figure}
\begin{center}
\begin{tabular}{ccc}
\begin{tikzpicture}
 [scale=.6,auto=left,every node/.style={circle,draw=black!50,fill=black!20}]
  \node (n1) at (0,0)  {1};
  \node (n2) at (2,2)  {2};
  \node (n3) at (2,0) {3};
    \node (n5) at (4,0)  {5};
  \node (n6) at (6,0)  {6};
   \node (n7) at (8,0)  {7};
  \node (n4) at (2,-2)  {4};

  \foreach \from/\to in {n1/n3, n2/n3,n3/n4,n3/n5,n5/n6,n6/n7}
    \draw (\from) -- (\to);
\end{tikzpicture}
& \qquad \qquad &
\begin{tikzpicture}
 [scale=.6,auto=left,every node/.style={circle,draw=black!50,fill=black!20}]
   \node (n1) at (0,0)  {1};
  \node (n2) at (2,2)  {2};
  \node (n3) at (2,0) {3};
    \node (n5) at (4,0)  {5};
  \node (n6) at (6,0)  {6};
   \node (n7) at (2,4)  {7};
  \node (n4) at (2,-2)  {4};

  \foreach \from/\to in {n1/n3, n2/n3,n3/n4,n3/n5,n5/n6,n2/n7}
    \draw (\from) -- (\to);
\end{tikzpicture} \\
$G_1$ &  & $G_2$
\end{tabular}
 \caption{Two graphs with different mp-charts, but with equal means and variances.} \label{fig-eqmv}
\end{center}
\end{figure}
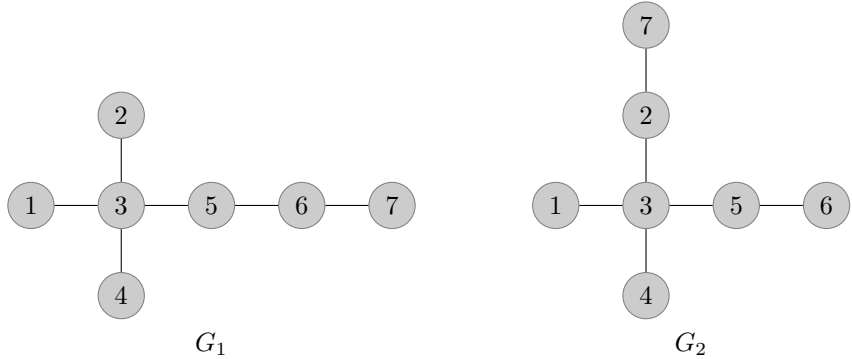
%
%The first graph $G_1$ has the edge set
%\begin{equation*}
%E_1 = \{(1,3), (2,3), (3,4), (3,5), (5,6), (6,7)\} \, ,
%\end{equation*}
%while the second graph $G_2$ has the edge set
%\begin{equation*}
%E_2=\{(1,3), (2,3), (3,4), (3,5), (5,6), (2,7)\} \, .
%\end{equation*}
Notice that $G_1$ and $G_2$ differ by only one edge. The two
mp-charts $H_{G_1}$ and $H_{G_2}$ have the same mean and variance,
namely
\begin{equation*}
{\mathbb E}(H_{G_1})={\mathbb E}(H_{G_2}) = \frac {12} {7} \qquad
\qquad {\mathbb V}(H_{G_1})={\mathbb V}(H_{G_2}) = \frac {170}
{147}
\end{equation*}
but the mp-charts are different:
\begin{eqnarray*}
H_{G_1}&=&(678,1512,1716,840,294,0,0,0) \\
H_{G_2}&=&(674,1480,1792,840,218,32,4,0) \, .
\end{eqnarray*}
\end{example}

The results above lead to a central limit theorem.
\begin{theorem} \label{clt}
Let $G_N$ be a graph with $N$ vertices and let $A_{G_N}$ be its
adjacency matrix. Let $\pi$ be a random permutation of $\{1 ,
\ldots , N\}$ chosen with uniform probability and define
\begin{equation}
S_N(\pi) = \sum_{i=1}^N (A_{G_N})_{i,\pi(i)} \, .
\end{equation}
If $\mathbb{V}(S_N)$ goes to infinity as $N \rightarrow \infty$,
then the distribution of $S_N$ is asymptotically normal.
\end{theorem}
\begin{proof}
We make use of Theorem 3 in \cite{hoeffding:51}. Define the
auxiliary matrix $R$ with elements
\begin{equation}
R_{i,j} = (A_{G_N})_{i,j} - \frac {d_i} N - \frac {d_j} N + \frac
1 {N^2} \sum_{h,k} (A_{G_N})_{h,k}
\end{equation}
Then, a sufficient condition for the asymptotic normality is that
\begin{equation} \label{suffcond}
\lim_{N \rightarrow \infty} \frac {\max_{1 \leq i,j \leq N}
R^2_{i,j} } { \frac 1 N \sum_{i,j=1}^N R^2_{i,j} } = 0 \, .
\end{equation}

Now observe that the numerator is bounded, as $-2 \leq R_{i,j}
\leq 2$ for all $i$ and $j$. Moreover, Theorem 2 in the same paper
\cite{hoeffding:51} states that
\begin{equation}
\mathbb{V}(S_N) = \frac 1 {N-1} \sum_{i,j=1}^N R_{i,j}^2 \, .
\end{equation}
Combining these facts, the result follows.
\end{proof}

\begin{remark}
Note that in our problem one can not use the classical central
limit theorems based on $\alpha$-mixing sequences or $m$-dependent
variables, see for instance \cite[Ch. 27]{billingsley:95} and
\cite{serfling:68}. Indeed, the covariance between $A_{i,\pi(i)}$
and $A_{j,\pi(j)}$ does not vanish as $|i-j|$ goes to infinity.
\end{remark}

In order to inspect the behavior of the convergence to the
Gaussian distribution, we have computed the mp-chart for some
graphs with $20$ vertices.

\begin{figure}
\begin{center}
\begin{tabular}{c}
\begin{tikzpicture}
  [scale=.45,auto=left,every node/.style={circle,draw=black!50,fill=black!20}]
  \node (n1) at (12,0) {};
  \node (n2) at (0,3)  {};
  \node (n3) at (3,3)  {};
  \node (n4) at (6,3)  {};
  \node (n5) at (9,3)  {};
  \node (n6) at (12,3)  {};
  \node (n7) at (15,3)  {};
  \node (n8) at (18,3)  {};
  \node (n9) at (21,3)  {};
  \node (n10) at (24,3) {};
  \node (n11) at (0,6)  {};
  \node (n12) at (3,6)  {};
  \node (n13) at (6,6)  {};
  \node (n14) at (9,6)  {};
  \node (n15) at (12,6)  {};
  \node (n16) at (15,6)  {};
  \node (n17) at (18,6)  {};
  \node (n18) at (21,6)  {};
  \node (n19) at (24,6)  {};
  \node (n20) at (12,9)  {};

  \foreach \from/\to in {n1/n2,n1/n3,n1/n4,n1/n5,n1/n6,n1/n7,n1/n8,n1/n9,n1/n10,
                         n11/n20,n12/n20,n13/n20,n14/n20,n15/n20,n16/n20,n17/n20,n18/n20,n19/n20,
                         n11/n2,n12/n3,n13/n4,n14/n5,n15/n6,n16/n7,n17/n8,n18/n9,n19/n10}
    \draw (\from) -- (\to);

\end{tikzpicture} \\
\begin{tabular}{cc}
\epsfig{file=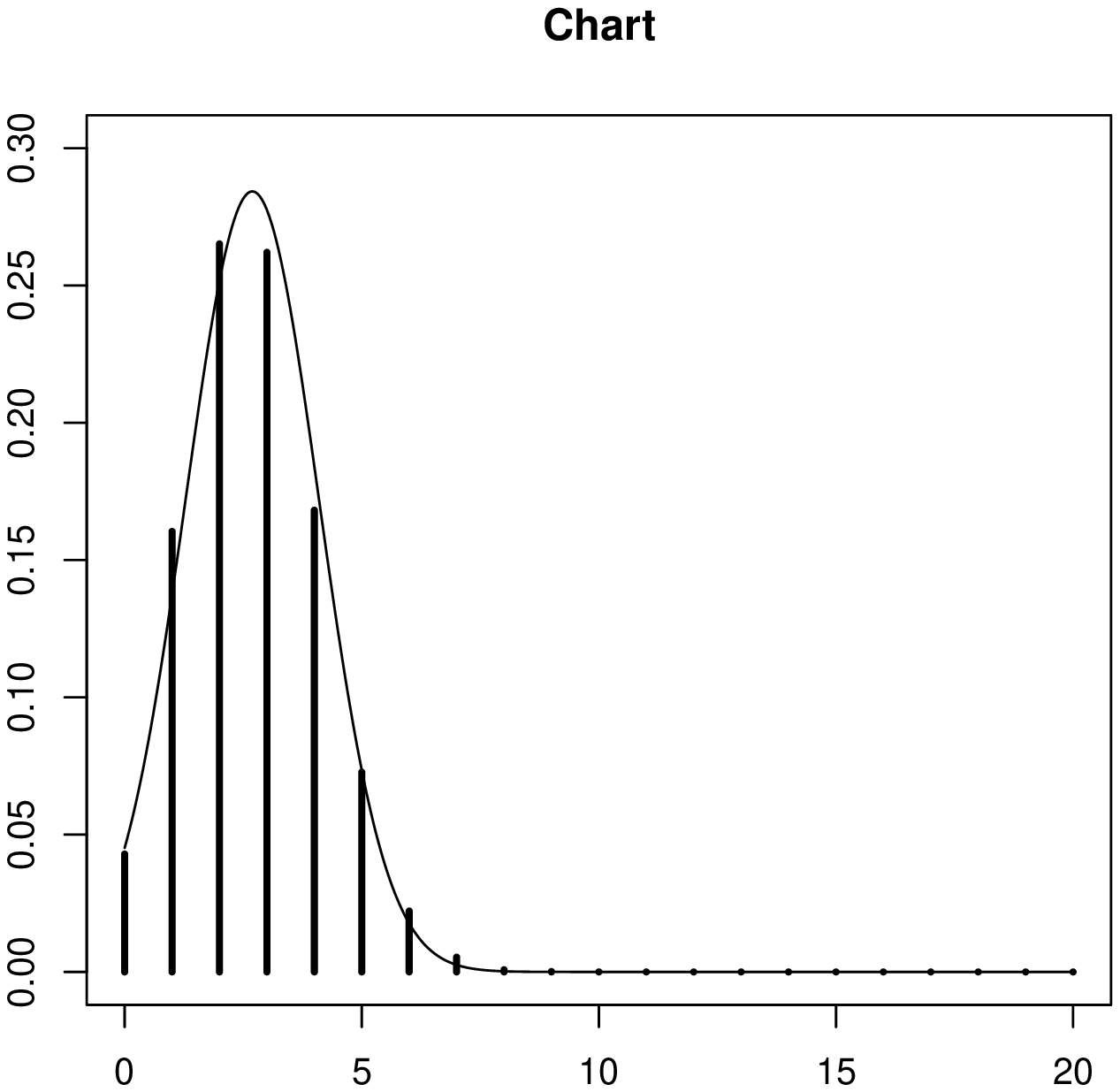,height=6cm,width=6cm}   &
\epsfig{file=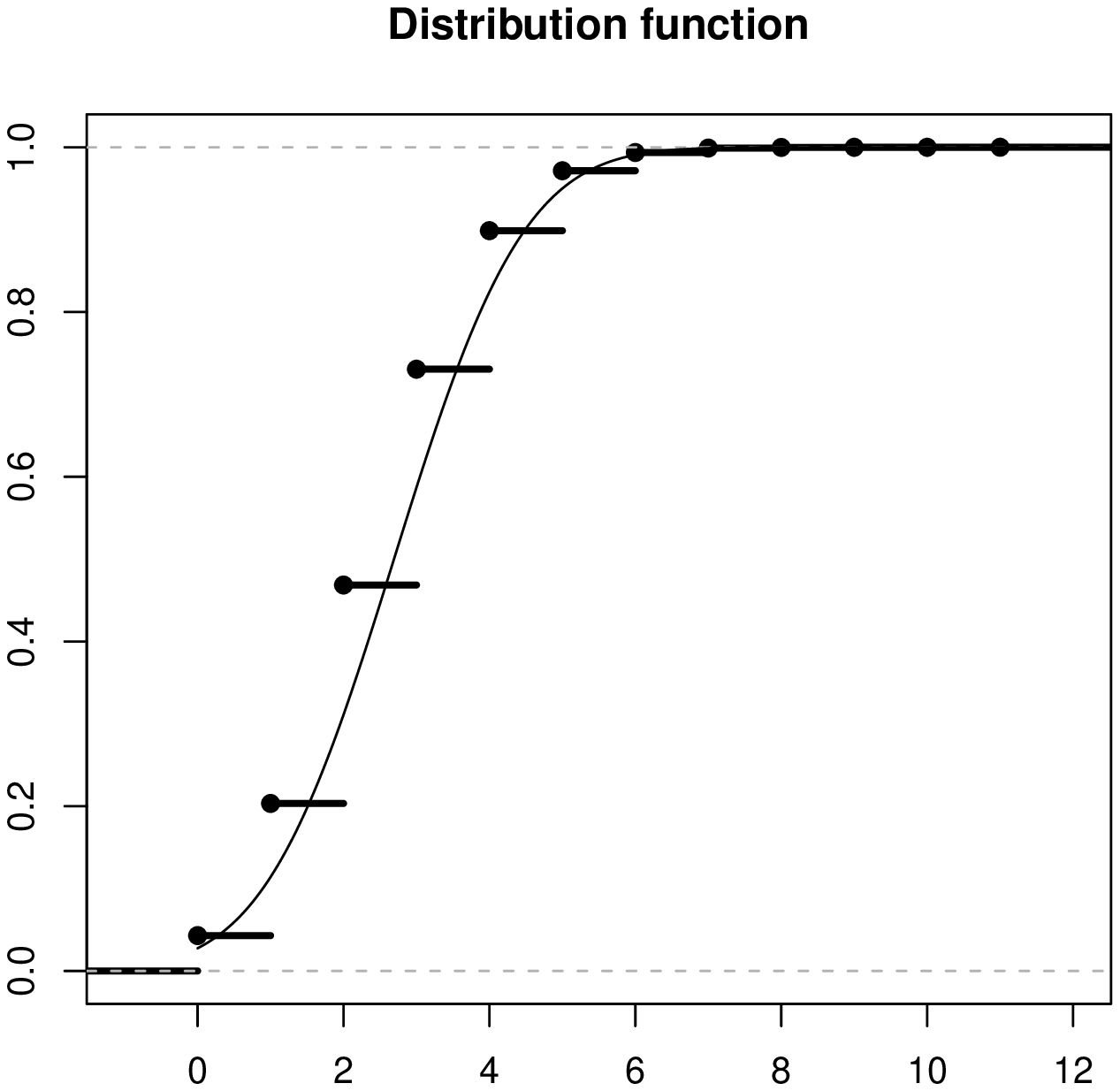,height=6cm,width=6cm}
\end{tabular}
\end{tabular}
\caption{Checking the convergence for $N=20$: First graph and its
results.}  \label{conv1}

\end{center}
\end{figure}

\begin{figure}
\begin{center}
\begin{tabular}{c}
\begin{tikzpicture}
  [scale=.6,auto=left,every node/.style={circle,draw=black!50,fill=black!20}]
  \node (n1) at (0:5cm) {};
  \node (n2) at (18:5cm)  {};
  \node (n3) at (36:5cm)  {};
  \node (n4) at (54:5cm)  {};
  \node (n5) at (72:5cm)  {};
  \node (n6) at (90:5cm)  {};
  \node (n7) at (108:5cm)  {};
  \node (n8) at (126:5cm)  {};
  \node (n9) at (144:5cm)  {};
  \node (n10) at (162:5cm) {};
  \node (n11) at (180:5cm)  {};
  \node (n12) at (198:5cm)  {};
  \node (n13) at (216:5cm)  {};
  \node (n14) at (234:5cm)  {};
  \node (n15) at (252:5cm)  {};
  \node (n16) at (270:5cm)  {};
  \node (n17) at (288:5cm)  {};
  \node (n18) at (306:5cm)  {};
  \node (n19) at (324:5cm)  {};
  \node (n20) at (342:5cm)  {};

  \foreach \from/\to in
  {n1/n2,n2/n3,n3/n4,n4/n5,n5/n6,n6/n7,n7/n8,n8/n9,n9/n10,n10/n11,
  n11/n12,n12/n13,n13/n14,n14/n15,n15/n16,n16/n17,n17/n18,n18/n19,n19/n20,n20/n1}
    \draw (\from) -- (\to);

  \foreach \from/\to in
  {n1/n3,n3/n5,n5/n7,n7/n9,n9/n11,n11/n13,n13/n15,n15/n17,n17/n19,n19/n1}
   \draw (\from) to [bend left=45] (\to);

  \foreach \from/\to in
  {n2/n4,n4/n6,n6/n8,n8/n10,n10/n12,n12/n14,n14/n16,n16/n18,n18/n20,n20/n2}
   \draw (\from) to [bend right=45] (\to);

\end{tikzpicture}   \\
\begin{tabular}{cc}
\epsfig{file=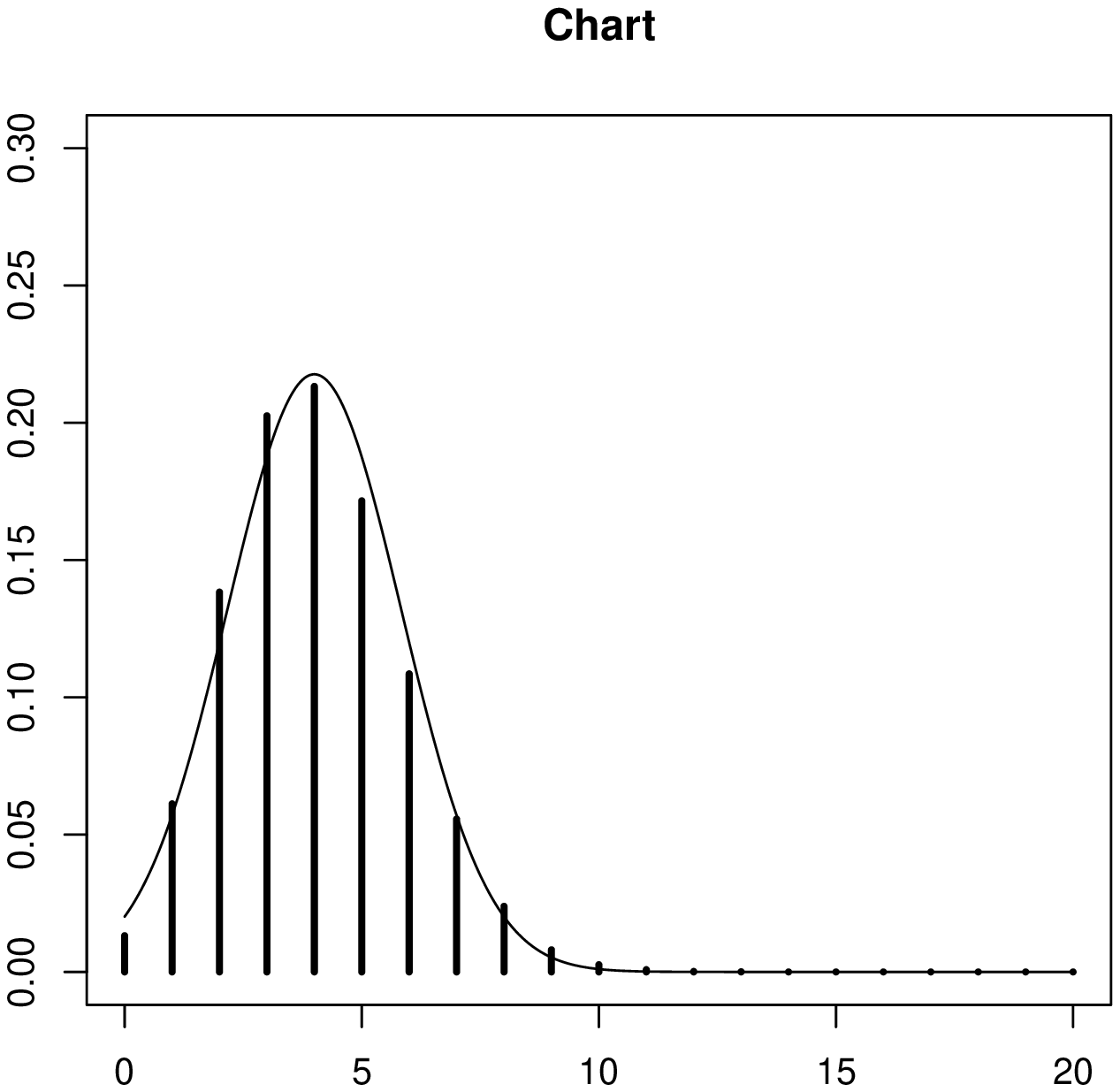,height=6cm,width=6cm}   &
\epsfig{file=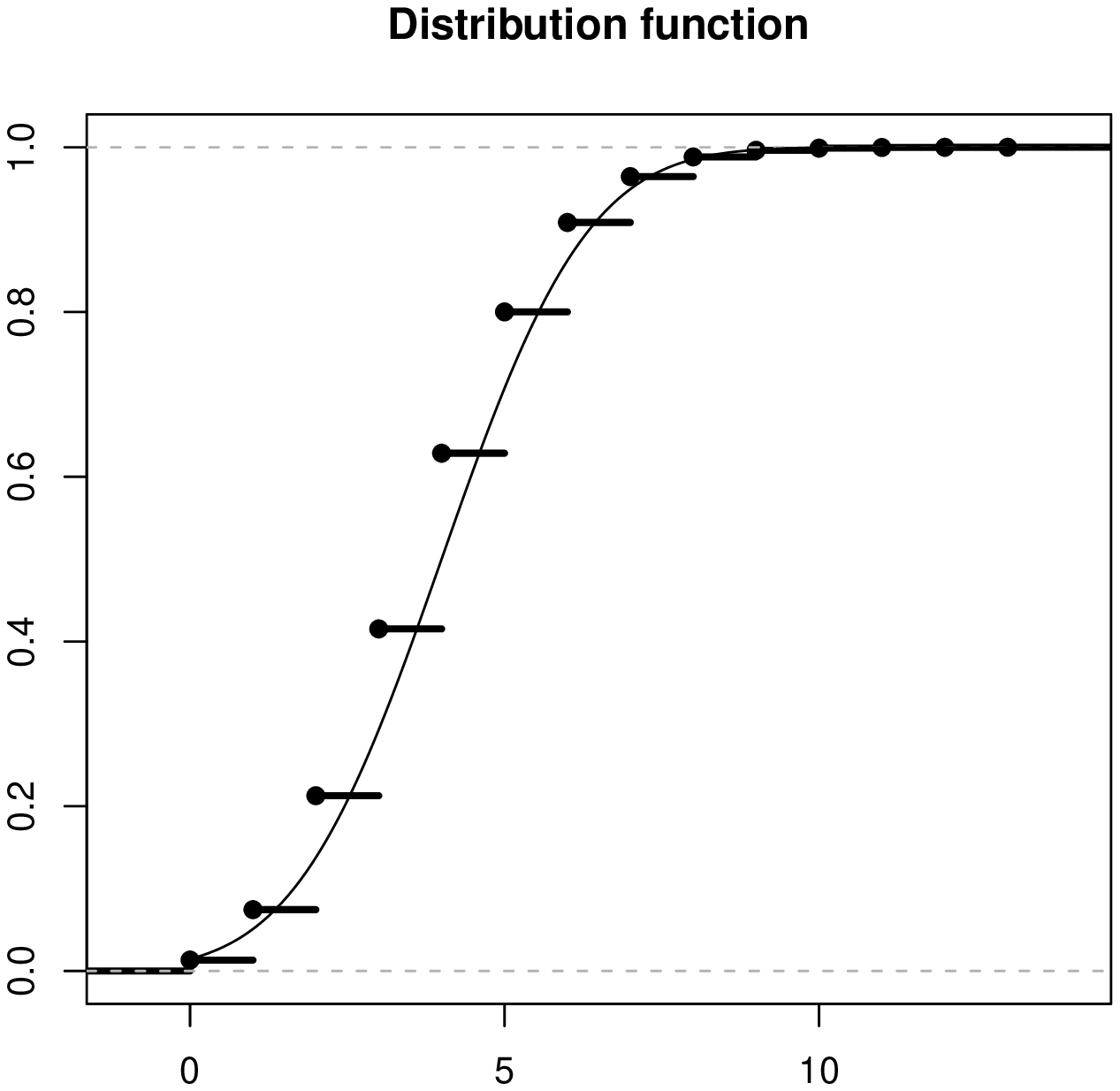,height=6cm,width=6cm}
\end{tabular}
\end{tabular}
\caption{Checking the convergence for $N=20$: Second graph and its
results.}  \label{conv2}
\end{center}
\end{figure}

\begin{figure}
\begin{center}
\begin{tabular}{cc}
\epsfig{file=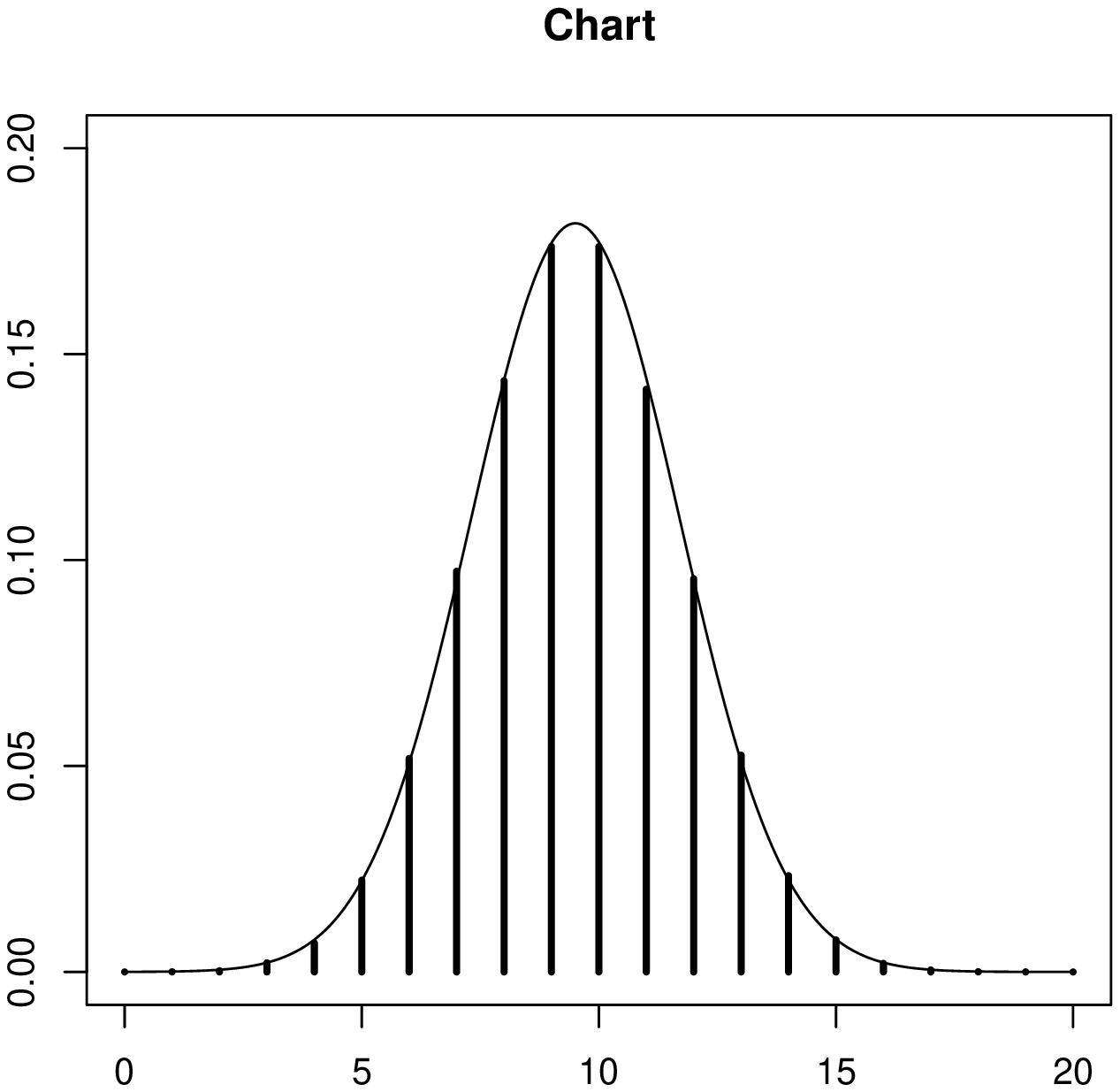,height=6cm,width=6cm}   &
\epsfig{file=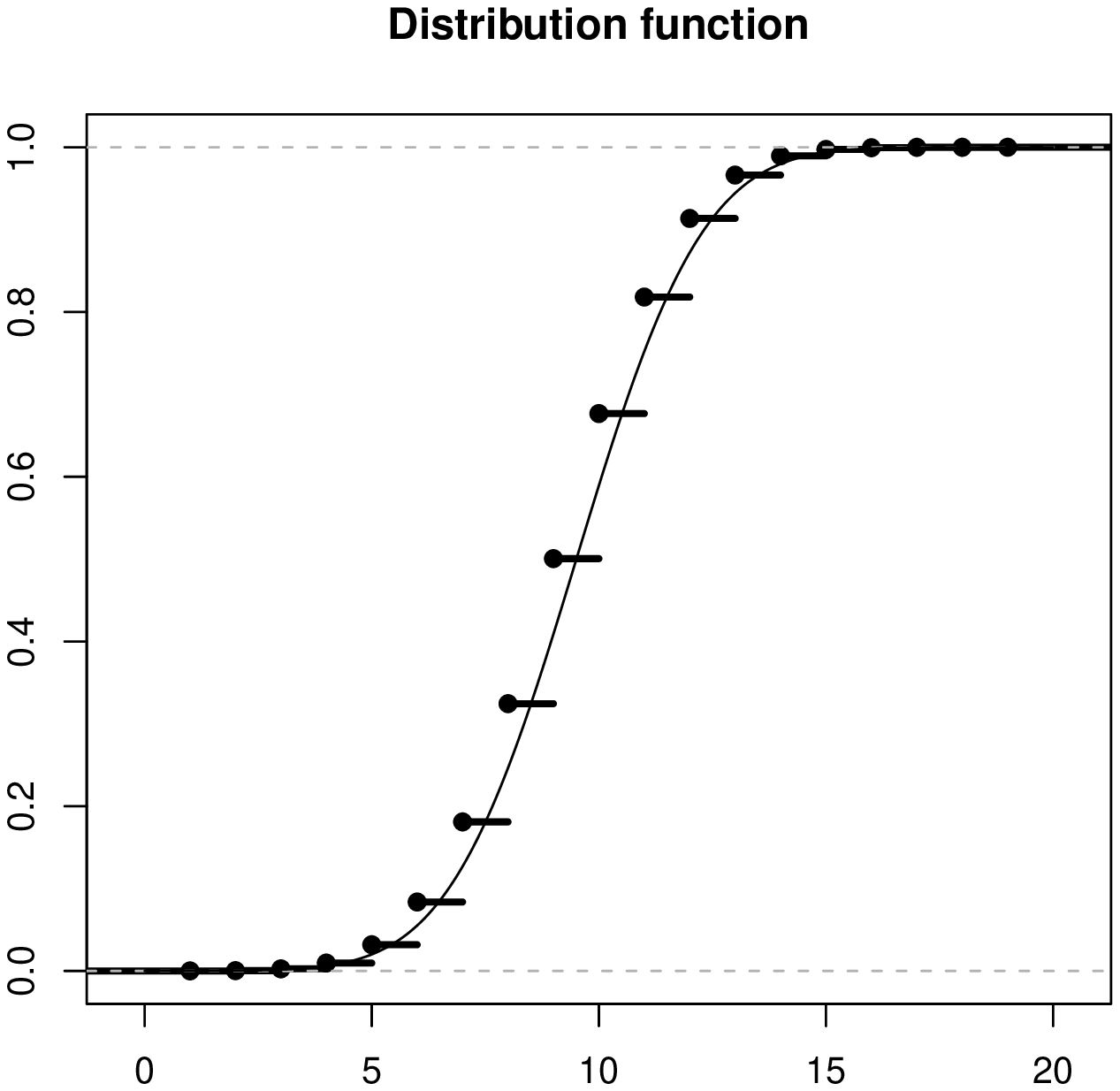,height=6cm,width=6cm}
\end{tabular}
\caption{Checking the convergence for $N=20$: The results for the
third graph.} \label{conv3}
\end{center}
\end{figure}

The three examples in Figures \ref{conv1}-\ref{conv3} show that
the convergence is quite good, meaning that the Gaussian
approximation is valid also for medium-sized graphs.

Two remarks are needed to understand the examples: $(a)$ The
mp-chart is approximated through a standard Monte Carlo technique,
sampling $100,000$ random permutations. This number is
considerably smaller than the total number of permutations ($20!
\cong 10^{18}$), but it provides quite accurate approximations;
$(b)$ The results are presented through two plots, showing the
mp-chart (normalized to $1$) and its distribution function, both
compared with the appropriate Normal distribution.

The first graph corresponds to an adjacency matrix with block
structure. The graph and the two plots of the results are
presented in Figure \ref{conv1}.

The second graph has a different shape, as it corresponds to an
adjacency matrix with band structure. The graph and the two plots
of the results are presented in Figure \ref{conv2}.

For the third graph analyzed here, we present only the results.
The graph has been constructed with $N(N-1)/4=95$ edges randomly
chosen among the $190$ edges of the complete graph with uniform
probability. The results are shown in Figure \ref{conv3}.

\section{The mp-chart of the complement graph}  \label{complsect}

In the literature, the complement of a graph $G=(V,E)$ is a graph
on the same vertex set $V$ and the set of edges $V^2 \setminus E$.
Since our starting graph has no loop $(e,e)$, this forces $G^c$ to
contain all of them. To avoid this problem we give a different
definition of complement graph, more useful for our purposes.

\begin{definition}
Given a graph $G=(V,E)$, its complement graph $G^c$ is  a graph on
the same vertex set $V$ and the set of edges $(V^2\setminus
\Delta) \setminus  E$, where $\Delta \subset E^2$ is the diagonal
set, i.e. $\Delta=\{(v,v) : v \in V \}$.
\end{definition}

\begin{remark}\label{latigrafi}
From the previous definition, we notice that, if $G$ and $G^c$
have respectively $\ell_G$ and $\ell_{G^c}$ edges, then
$\ell_G+\ell_{G^c}={N(N-1)}/{2}$.
\end{remark}

As mentioned in Section \ref{mpchartsect}, there are some nice
properties linking the mp-chart $H_G$ of a graph $G$ with the
mp-chart of its complement $G^c$. To study these connections,  we
start  with a preliminary lemma.

\begin{lemma}\label{LemTeTc}
Let $G$ be a graph with $N$ vertices. Denote by $T_{ij}^c$ the
scalar product between the $i$-th row and $j$-th column in the
adjacency matrix of $G^c$. The following formula relates the
quantities $T_{ij}^c$ and $T_{ij}$:
\begin{equation}\label{TeTc}
 T^c_{ij}=T_{ij}+N-2-d_i-d_j+<r_i,E_j>+<r_j,E_i>.
\end{equation}
where the $E_i$'s are the vectors in the canonical basis of
${\mathbb R}^N$.
\end{lemma}

\begin{proof}
Define the vector
\begin{equation}\label{vij}
v_{ij}=r_i-r_j-<r_i-r_j,E_i>E_i-<r_i-r_j,E_j>E_j
\end{equation}
Since $<r_i-r_j,E_i>=<r_i,E_i>-<r_j,E_i>$ and $<r_t,E_t>=0$ (in
fact, the $t$-th coordinate of $r_t$ is zero, while $E_t$ has a
$1$ in the $t$-th coordinate and $0$ elsewhere), we can write
\begin{equation*}
v_{ij}=r_i-r_j+<r_j,E_i>E_i-<r_i,E_j>E_j
\end{equation*}
The scalar product $<v_{ij},v_{ij}>$ measures the number of
positions, out of the diagonal, where $r_i$ and $r_j$ are
different. Thus, if we denote by $T^c_{ij}$ the scalar product of
the corresponding lines, $r_i^c$, $r_j^c$ in the complement graph,
one has
\begin{equation}\label{passaggio}
T^c_{ij}=N-2-T_{ij}-<v_{ij},v_{ij}>.
\end{equation}
Now, we substitute in the previous formula the expression of
$v_{ij}$ given in (\ref{vij}), and we obtain:
\begin{equation*}
\begin{split}
T^c_{ij}&= N-2-T_{ij}-<v_{ij},v_{ij}>  =  N-2-T_{ij}-\\
&-<r_i-r_j+<r_j,E_i>E_i-<r_i,E_j>E_j,r_i-r_j+<r_j,E_i>E_i-<r_i,E_j>E_j>.
% &=N-2-T_{ij}-<r_i,r_i-r_j+<r_j,E_i>E_i-<r_i,E_j>E_j>+\\
% &+<r_j,r_i-r_j+<r_j,E_i>E_i-<r_i,E_j>E_j>+\\
% &-<r_j,E_i>E_i,r_i-r_j+<r_j,E_i>E_i-<r_i,E_j>E_j>+\\
% &+<r_i,E_j>E_j,r_i-r_j+<r_j,E_i>E_i-<r_i,E_j>E_j=\\
% &=N-2-T_{ij}-<r_i,r_i>+<r_i,r_j>-<r_i,<r_j,E_i>E_i>+<r_i,<r_i,E_j>E_j>+\\
% &+<r_j,r_i>-<r_j,r_j>+<r_j,<r_j,E_i>E_i>-<r_j,<r_i,E_j>E_j>+\\
% &-<<r_j,E_i>E_i,r_i>+<<r_j,E_i>E_i,r_j>-<<r_j,E_i>E_i,<r_j,E_i>E_i>+\\
%&+ <<r_j,E_i>E_i,<r_i,E_j>E_j>+<<r_i,E_j>E_j,r_i>-<<r_i,E_j>E_j,r_j>+\\
%&+<<r_i,E_j>E_j,<r_j,E_i>E_i>-<<r_i,E_j>E_j,<r_i,E_j>E_j>=
\end{split}
\end{equation*}
Noting that $<r_t,r_t>=d_t$ and $<r_i,r_j>=<r_j,r_i>=T_{ij}$,  a
straightforward computation leads to
\begin{equation*}
%\begin{split}
T^c_{ij}%&=N-2+T_{ij}-d_i+<r_i,<r_i,E_j>E_j>-<r_j,r_j>+<r_j,<r_j,E_i>E_i>+\\
 %&+<<r_j,E_i>E_i,r_j>-<<r_j,E_i>E_i,<r_j,E_i>E_i>+\\
 %&+<<r_i,E_j>E_j,r_i>-<<r_i,E_j>E_j,<r_i,E_j>E_j>=\\
 %&=N-2+T_{ij}-d_i-d_j+<r_i,E_j>^2+<r_j,E_i>^2+\\
 %&+<r_j,E_i>^2-<r_j,E_i>^2+<r_i,E_j>^2-<<r_i,E_j>^2=\\
%&
= T_{ij}+N-2-d_i-d_j+<r_i,E_j>^2+<r_j,E_i>^2.
% \end{split}
\end{equation*}
Since, for all $i$,$j$, with $i \not=j$,  the value of $<r_i,E_j>$
can be either $0$ or $1$, we can remove the squares from the
previous formula, leading to Equation (\ref{TeTc}).
\end{proof}

\begin{theorem}\label{compl}
Given a graph $G$ with $N$ vertices,
\begin{itemize}
\item[(a)] $\mathbb{E}(H_{G^c})=N-1-\mathbb{E}(H_G)$;

\item[(b)]
$\mathbb{V}(H_{G^c})=\mathbb{V}(H_G)+1-2{\mathbb{E}(H_G)}/{(N-1)}$.
\end{itemize}
\end{theorem}

\begin{proof}
To prove part $(a)$, it is enough to use Theorem \ref{Media} and
Remark \ref{latigrafi}. One has
\begin{equation*}
\mathbb{E}(H_{G^c})=\frac{2\ell_{G^c}}{N}=\frac{N(N-1)-2\ell_G}{N}=N-1-\mathbb{E}(H_G)
\, .
\end{equation*}

To prove part $(b)$, we apply Theorem \ref{Varianza} to the graph
$G^c$. Therefore we have:
\begin{equation*}
\mathbb{V}(H_{G^c})= \sum_{i=1}^N \frac {d_i^c(N-d_i^c)} {N^2} +
\sum_{\substack{i,j=1 \\ i \ne j}}^N \frac {d_i^cd_j^c -
NT^c_{i,j}} {N^2(N-1)} \, .
\end{equation*}
The degree $d_i^c$ of a vertex in the complement graph is given by
$d_i^c=N-1-d_i$. Using also Lemma \ref{LemTeTc} one can write
$\mathbb{V}(H_{G^c})$ as

\begin{equation} \label{lunga}
\begin{split}
&\mathbb{V}(H_{G^c}) = \sum_{i=1}^N \frac {(N-1-d_i)(d_i+1)} {N^2} + \\
&+\sum_{ \substack{i,j=1 \\ i \ne j}}^N \frac {(N-1-d_i)(N-1-d_j)-
N(T_{ij}+N-2-d_i-d_j+<r_i,E_j>+<r_j,E_i>)} {N^2(N-1)}
\end{split}
\end{equation}

The first sum in formula \eqref{lunga} can be written as
\begin{equation*}
\sum_{i=1}^N \frac {(d_i)(N-d_i)} {N^2} + \sum_{i=1}^N \frac
{(N-1-2d_i)} {N^2} \, ,
\end{equation*}
where an easy computation shows that
\begin{equation*}
\sum_{i=1}^N \frac {(N-1-2d_i)} {N^2} =\frac
{N-1}{N}-\frac{4\ell_G}{N^2} \, .
\end{equation*}

About the second sum in formula \eqref{lunga}, we observe that
each term can be expressed as

\begin{equation*}
\begin{split}
&\frac {(N-1-d_i)(N-1-d_j)- N(T_{ij}+N-2-d_i-d_j+<r_i,E_j>+<r_j,E_i>)} {N^2(N-1)} =\\
&= \frac {d_id_j -
NT_{i,j}+d_i+d_j+1-N(<r_i,E_j>+<r_j,E_i>)}{N^2(N-1)} \, .
\end{split}
\end{equation*}

Since
\begin{equation*}
\sum_{\substack{i,j=1 \\ i \ne j}}^N <r_i,E_j>+<r_j,E_i>) =4\ell_G
\ \mbox{ and } \ \sum_{\substack{i,j=1 \\ i \ne j}}^N
(d_i+d_j)=4\ell_G (N-1)
\end{equation*}
the second sum becomes
\begin{equation*}
 \sum_{\substack{i,j=1 \\ i \ne j}}^N \frac {d_id_j - NT_{i,j}} {N^2(N-1)}
+\frac{4(N-1)\ell_G-N(N-1)-4N\ell_G}{N^2(N-1)} \, .
\end{equation*}
Thus, we obtain
$$\mathbb{V}(H_{G^c})=\mathbb{V}(H_G)+1-\frac{4\ell_G}{N(N-1)}$$

\noindent and, considering Theorem \ref{Media}, the formula in (b)
follows.
\end{proof}

Some major remarks on the Theorem above are now in order.

\begin{remark}
As a first trivial example, we consider the limit situation of an
empty graph. Let $G$ be the empty graph. Its mp-chart has
$\mathbb{E}(H_G)=\mathbb{V}(H_{G})=0$. In this case $G^c$ is the
complete graph with ${N(N-1)}/{2}$ edges and, by Theorem
\ref{compl} part (b), one has $\mathbb{V}(H_{G^c})=1$. This can be
verified also by direct computation. As a matter of fact, the
adjacency matrix of $G^c$ consists of non-zero entries out of the
diagonal. Thus $d_i^c=N-1$ for all $i$ and $T_{ij}^c=N-2$ for all
$i$,$j$, with $i \not=j$. Hence
\begin{equation*}
\mathbb{V}(H_G)= \sum_{i=1}^N \frac {d_i(N-d_i)} {N^2} +
\sum_{\substack{i,j=1 \\ i \ne j}}^N \frac {d_id_j - NT_{i,j}}
{N^2(N-1)} =\frac{N-1}{N}+\frac{(N-1)^2-N(N-2)}{N^2(N-1)}=1 \, .
\end{equation*}
These computations show that the variance of the mp-chart of a
complete graph is invariant on the number of vertices.
\end{remark}

%\begin{Rmk}
%Let $G$ be a graph with $N$ vertices and $\ell_G$ edges. If we
%apply Theorem \ref{compl} to $G^c$ to compute the variance of $G$
%with expect to find the same result as to compute it directly. In
%fact, using part (b) of Theorem \ref{compl} both on $G$ and $G^c$
%one has
%\begin{equation*}
%\mathbb{V}(H_{G})=
%\mathbb{V}(H_{G^c})+1-\frac{2\mathbb{E}(H_G)}{N-1}=
%\mathbb{V}(H_G)+1-\frac{2\mathbb{E}(H_G)}{N-1}+1-\frac{2\mathbb{E}(H_{G^c})}{N-1}
%\end{equation*}
%and, by part (a) of Theorem \ref{compl}, we know that
%$2-\frac{2\mathbb{E}(H_G)}{N-1}-\frac{2\mathbb{E}(H_{G^c})}{N-1}=0$.
%\end{Rmk}

\begin{remark}
Few straightforward algebraic calculations show that the
difference ${\mathbb V}(H_{G^c})-{\mathbb V}(H_{G})$ lies between
$-1$ and $1$. Therefore, under the hypotheses of Theorem
\ref{clt}, when $N$ goes to infinity we have that ${\mathbb
V}(H_{G^c}) \cong {\mathbb V}(H_{G})$. Intuitively, the difference
between ${\mathbb V}(H_{G^c})$ and ${\mathbb V}(H_{G})$ depends on
the diagonal entries which are forced to be zero in the adjacency
matrix $A_{G^c}$. The effect of these entries vanishes when the
size of the graph goes to infinity.
\end{remark}

\begin{remark}
Another interesting property follows from Theorem \ref{compl}.
First, notice that ${\mathbb E}(H_{G^c})={\mathbb E}(H_{G})$
implies that $G$ and $G^c$ have the same number of edges, namely
$N(N-1)/4$. (This is not possible for all values of $N$). In such
a case, $H_{G^c}$ and $H_G$ are forced to have the same variance,
no matter how is complicated the graph $G$.
\end{remark}

The computation of the whole mp-chart of the complement graph
$G^c$ from the mp-graph of $G$ is less easy. Given a graph G, we
build a $(N+1)\times (N+1)-$matrix $M_G$, indexed, both on rows
and columns, by $\{0, \dots, N \}$, an defined as follows. The
entry $(M_G)_{i,j}$ is the number  of permutations $\pi$ such that
$\sum_{s=1}^N (A_G)_{s,\pi(s)}=j$ and $\pi$ has $i$ diagonal
elements (that is, $\pi(s)=s$ for $i$ elements).

\begin{example}\label{matrixgraph} Consider the graph $G$ with matrix
\begin{equation*}
A_G= \begin{pmatrix}
0 & 0 & 0 &1\\
0 & 0 & 1 &0\\
0& 1&0 &1 \\
1& 0& 1 &0
\end{pmatrix} \, .
\end{equation*}
The corresponding matrix $M_G$ is
\begin{equation*}
M_G=  \begin{pmatrix}
1 & 2 & 3 &2&1\\
0 & 4 & 4 & 0 &0\\
3& 0&3 &0 &0\\
0& 0& 0 &0 &0 \\
1 &0 &0 &0 &0
\end{pmatrix} \, .
\end{equation*}
\end{example}

The matrix $M_G$ allows the computation of both the mp-charts
$H_G$ and $H_{G^c}$. Roughly speaking, to compute the mp-chart of
$G$ is is enough to sum the columns of $M_G$, while to compute the
mp-chart of $G^c$ we need to sum the entries of suitable diagonals
of $M_G$. More precisely, the following relations hold true.

\begin{proposition} \label{prop-complem}
For a graph $G$, we have for all $j=1, \ldots,N$:
\begin{itemize}
\item[(a)] The components of the mp-chart $H_G$ are:
\begin{equation*}
h_G(j)=\sum_{i=0}^N (M_G)_{i,j} \, ;
\end{equation*}
\item[(b)] The components of the mp-chart $H_{G^c}$ are:
\begin{equation*}
h_{G^c}(j)=\sum_{i=0}^{N-j} (M_G)_{i,N-j-i}.
\end{equation*}
\end{itemize}
\end{proposition}
\begin{proof}
The first relation follows by the definition of mp-chart, as the
sum of entries in the $j$-th column of $M_G$ is the number of
permutations with $j$ elements equal to $1$, that is $h_G(j)$.

To prove the second relation, it is enough to prove that for all
$i$ and $j$ we have: $(M_{G^c})_{i,j}=(M_G)_{i,N-j-i}$. Suppose
that $\pi$ is such that $\sum_{s=1}^N (A_G)_{s,\pi(s)}=j$ and
$\pi$ has $i$ diagonal elements. When we consider $\pi$ on
$A(G^c)$, we have $A(G^c)_{s,\pi(s)}=1$ for the $s$ such that
$(A_G)_{s,\pi(s)}=0$, except for the diagonal entries where we
still have 0. Hence there are $N-j-i$ entries in $A(G^c)$ such
that $A(G^c)_{s,\pi(s)}=1$.  This completes the proof.
\end{proof}

\begin{example}
Consider the complement graph $G^c$ of the graph $G$ in Example
\ref{matrixgraph}:
\begin{equation*}
A_{G^c}= \begin{pmatrix}
0 & 1 & 1 &0\\
1 & 0 & 0 &1\\
1& 0&0 &0 \\
0& 1& 0 &0
\end{pmatrix} \, .
\end{equation*}
Then
\begin{equation*}
M_{G^c}=  \begin{pmatrix}
1 & 2 & 3 &2&1\\
0 & 4 & 4 & 0 &0\\
3& 0&3 &0 &0\\
0& 0& 0 &0 &0 \\
1 &0 &0 &0 &0
\end{pmatrix} \, .
\end{equation*}
We notice that in this simple case $M_G=M_{G^c}$. This is due to
the fact that $G$ and $G^c$, up to the labels of the vertices, are
equivalent.
\end{example}

\begin{remark}
In principle, the max-plus permanent of the complement graph $G^c$
can be computed from the matrix $M_G$. In fact, Proposition
\ref{prop-complem} shows that the mp-chart of the complement graph
can be computed from the matrix $M_G$ by adding along suitable
diagonals and the max-plus permanent is just the position of the
last non-zero element of the mp-chart. However, as the matrix
$M_G$ is not easy to compute for large graphs, this approach does
not help in actual computations.
\end{remark}

\section{Future directions}  \label{futuresect}

The max-plus permanent and the mp-chart studied in this paper lead
to several new questions. In fact, we have analyzed here only the
case of undirected graph. Therefore, among the future directions
of our research, there will be the extension of the definition of
max-plus permanent and mp-chart for undirected and weighted
graphs. Moreover, special classes of graphs may be studied, such
as bipartite graphs, or fixed-degree graphs. In particular,
fixed-degree graphs correspond to adjacency matrices with fixed
margins and in that context algebraic and combinatorial methods
have demonstrated already their potential.

Another possible research direction is strictly in graph theory.
As a matter of fact it could be interesting to compare the
mp-chart of a graph with other well-known descriptors of its
complexity. In a recent work in progress the mp-chart is compared
to matching polynomials. By several examples we know that there
exist different graphs with the same mp-chart. We notice that, in
all these cases, also the matching polynomials coincide. So a
principal question would be: If two non-isomorphic graphs have the
same mp-chart, then their matching polynomials are equal?

Finally, applications to large graphs, possibly through simulation
techniques, will be investigated in order to use the tools
presented in this paper to real data examples.

%
%\section*{Acknowledgments}
%To do
%
%\bibliographystyle{siam}
%\bibliography{tuttopm}

\end{document}